\documentclass[12pt,reqno]{article}

\usepackage[usenames]{color}
\usepackage{amssymb}
\usepackage{graphicx}
\usepackage{amscd}

\usepackage[colorlinks=true,
linkcolor=webgreen,
filecolor=webbrown,
citecolor=webgreen]{hyperref}

\definecolor{webgreen}{rgb}{0,.5,0}
\definecolor{webbrown}{rgb}{.6,0,0}

\usepackage{color}
\usepackage{fullpage}
\usepackage{float}

\usepackage{psfig}
\usepackage{graphics,amsmath,amssymb}
\usepackage{amsfonts}
\usepackage{latexsym}

\DeclareMathOperator{\ld}{ld}

\begin{document}


\begin{center}
\vskip 1cm{\LARGE\bf Generalized Number Derivatives}
\vskip 1cm
\large
Michael Stay\\
Department of Computer Science\\
University of Auckland\\
Private Bag 92019\\
Auckland 1020\\
New Zealand \\
\href{mailto:staym@clear.net.nz}{\tt staym@clear.net.nz} \\
\end{center}

\vskip .2 in
\begin{abstract}
We generalize the concept of a number derivative, and examine one
particular instance of a deformed number derivative for finite field
elements.  We find that the derivative is linear when the deformation
is a Frobenius map and go on to examine some of its basic properties.

\end{abstract}

\newtheorem{theorem}{Theorem}[section]
\newtheorem{proposition}{Proposition}[section]
\newtheorem{corollary}{Corollary}[section]
\newtheorem{lemma}{Lemma}[section]

\newcommand{\ctr}[1]{\begin{center} #1 \end{center}}
\newcommand{\oo}[0]{\infty}
\newcommand{\al}[1]{\begin{align} #1 \end{align}}
\newcommand{\alx}[1]{\begin{align*} #1 \end{align*}}
\newcommand{\bigparens}[1]{\mbox{\huge (} #1 \mbox{\huge )}}
\renewcommand{\thefootnote}{\fnsymbol{footnote}}


\section{Introduction}
The so-called ``number derivative'' seems to have been invented
independently at least three times \cite{KOW02,CI03,UA03}.  Here we
present a generalization of the number derivative that applies to
nearly anything one might reasonably call a number.  Afterwards, we
examine the case of a specific number derivative on finite fields and
some of its basic properties.

We generalize the concept of a number derivative to the following
algorithm; in order to illustrate each step, we will present the
corresponding step from the standard number derivative, denoted $N$,
and our number derivative, denoted $S$.

The notation we use below requires some care. Multiplication is denoted by a dot $[\cdot]$ or by concatenation of symbols: $x\cdot x^2=xx^2=x^3$.  The notation $x^n$ denotes a {\em function}: 
    $$x^n(y)=y^n,$$ 
so $x(y)=y$ is the identity function and $x^0(y)=1$ for all $y$.  Parentheses, when preceded by a function or operator, denote composition or application, respectively:
    $$x^2(x^n) = (x^n)^2 = x^{2n}.$$
Application is left associative:
    $$f(g)(h)=(f(g))(h),$$
and takes precedence over multiplication:
    $$gh(f)\ne g(f)\cdot h(f).$$

\begin{enumerate}
\item Choose a parameterized canonical form.  In the case of $N$, this
consists of representing each integer as a product of prime powers; the
parameters are the primes.  In the case of $S$, we choose a generator
$\theta$ of the finite field GF($p^k$) and express each finite field
element as $\theta^n$.  Here, the parameter is $\theta$.

\item Convert this canonical form into a function.  The algorithm $N$
takes each prime power $p_i^{k_i}$ to a function
$x_i^{k_i}(y_i)=y_i^{k_i}$.  The algorithm $S$ replaces $\theta^n$ with
the function $x^n(y)=y^n$.

\item Differentiate the function with respect to the parameters.  The
algorithm $N$ computes $D(f) = (\sum_{i}\frac{\partial}{\partial
x_i})(f)$.  The algorithm $S$ computes the $s$-derivative $D_s(f)$.

\item Evaluate the derivative at some function of the parameters,
typically the identity function.  The algorithm $N$ computes
$D(f)|_{y_i = p_i}$.  The algorithm $S$ computes $D_s(f)|_{y=\theta}$.
\end{enumerate}

\section{Exponential quantum calculus}

We begin with the operator $M_s(f)=f(x^s)$.  The $s$-differential is then $d_s=M_s-x$ and the $s$-derivative is
    $$D_s(f) = \frac{d_s (f)}{d_s (x)}.$$
The $s$-derivative of an element $x^n(\theta)$ is
    $$D_s(x^n)(\theta)=\frac{M_s(x^n)-x^n}{M_s(x)-x}(\theta)
              =\frac{x^{ns}-x^n}{x^s-x}(\theta)=([n]x^{n-1})(\theta),$$
where $$[n] \equiv \frac{x^{(s-1)n}-x^0}{x^{s-1}-x^0}.$$

The $s$-deformation has many similarities to the $q$-deformation that
results in the quantum calculus \cite{KC02}.  To get the
$s$-deformation from the $q$-deformation, one replaces the constant $q$
by the function $x^{s-1}$.  Since this is the same transformation we
chose to use in the second step of the algorithm $S$, both derivatives
give rise to the same number derivative.

Since the notation is somewhat simpler for the $q$-derivative, we will adopt it through most of the paper.  The operator $M_q = M_s$:
    $$M_s(f)=f(x^s)=f(x^{s-1}x)=f(qx)=M_q(f).$$  
The $q$-differential is $d_q=M_q-x$ and the $q$-derivative is
    $$D_q(f) = \frac{d_q(f)}{d_q(x)}.$$
The $q$-derivative of an element $x^n(\theta)$ is
    $$D_q(x^n)(\theta)=\frac{M_q(x^n)-x^n}{M_q(x)-x}(\theta)
              =\frac{q^nx^n-x^n}{qx-x}(\theta)=([n]x^{n-1})(\theta),$$
where $$[n] \equiv \frac{q^n-q^0}{q^1-q^0}.$$

Also, in the portions of the paper directly concerning the algorithm
$S$, we will usually omit the final application of the functions to
$\theta$.

\section{Identities}

For what functions $q=x^{s-1}$, if any, is this number derivative
linear?  Let $\theta^a+\theta^b=\theta^c$.  Then
    $$D_q(x^c)(\theta)=\frac{x^{sc}-x^c}{x^{s}-x}(\theta)
      =\frac{(\theta^a + \theta^b)^{s} -
      \theta^c}{\theta^{s}-\theta}.$$ 
On the other hand,
    $$(D_q(x^a)+D_q(x^b))(\theta)=\frac{x^{as}+x^{bs}-(x^a+x^b)}{x^s-x}(\theta)
=\frac{\theta^{as}+\theta^{bs}-\theta^c}{\theta^s-\theta},$$ so
we want the cross terms in the binomial $(\theta^a + \theta^b)^{s}$ to
be zero modulo $p$.  This only occurs when $s$ is a power of $p$, so
the derivative is linear if and only if $M_s$ is a Frobenius map.  In
the rest of the paper, we will only consider $q=x^{s-1}$ of this form.

The derivation of the product rule is the same as that for the $q$-derivative:
\al{
    D_q(fg)&= \frac{M_q(fg) - fg}{M_qx-x} \nonumber \\
                 &= \frac{M_q(f)M_q(g)-M_q(g)f+M_q(g)f-fg}
                         {M_q(x)-x} \nonumber \\
                 &= M_q(g)\frac{M_q(f)-f}{M_q(x)-x}+
                    f\frac{M_q(g)-g}{M_q(x)-x} \nonumber \\
                 &= M_q(g)D_q(f) + fD_q(g)   \label{eq:prod1} \\
                 &= gD_q(f) + M_q(f)D_q(g), \label{eq:prod2}
}
where (\ref{eq:prod2}) follows by symmetry.

The same is true for the quotient rule.  Since by (\ref{eq:prod1}),
\al{
    D_q (f) &= D_q(g\frac{f}{g}) \nonumber \\
             &= M_q(g)D_q(\frac{f}{g}) + \frac{f}{g}D_q(g), \nonumber 
    \intertext{we have}
    D_q (\frac{f}{g}) &=\frac{gD_q(f)-fD_q(g)}{gM_q(g)} \label{eq:quot1}\\
                    &=\frac{M_q(g)D_q(f)-M_q(f)D_q(g)}{gM_q(g)}\label{eq:quot2}
}
where (\ref{eq:quot2}) follows from (\ref{eq:prod2}) instead.

Note that while there is not a general chain rule for the standard
$q$-derivative, we can use the fact that every element is of the form
$x^n(\theta)$ to find one for this derivative:
\alx{
    D_q(g(x^n))&=\frac{M_q(g(x^n))-g(x^n)}{M_q(x)-x}\cdot \frac{M_q(x^n)-x^n}{M_q(x^n)-x^n} \\
               &=\frac{M_q(g(x^n))-g(x^n)}{M_q(x^n)-x^n}\cdot \frac{M_q(x^n)-x^n}{M_q(x)-x} \\
               &=\frac{M_q(g(x^n))-g(x^n)}{M_q(x(x^n))-x(x^n)}\cdot D_q(x^n) \\
               &=D_q(g)(x^n)\cdot D_q(x^n)                   
}
While the product and quotient rules (\ref{eq:prod1})-(\ref{eq:quot2})
are the same as those typically given \cite{KC02}, this rule differs:
since $q$ is the function $x^{p^j-1}$ instead of a constant, we
evaluate it at $x^n$ rather than take the $q^n$-derivative of $g$ in
the first term.

Finally, the $q$-numbers $[n]$ satisfy 
    $$[n+1] = q^0+q[n] \mbox{\;\;and\;\;} [n+1]-[n]=q^n.$$

\section{Constants}

Under what conditions does \al{d_q(x^n) &= 0 \label{eq:homogen}?}  We
have
    $$M_q(x^n) - x^n=0$$ 
    which implies
    $$q^nx^n = x^n$$ 
    and
    $$q^n = x^{n(p^j-1)} = x^0 = 1$$ if $n\ne -\oo$.  Therefore, $D_q
x^n = 0$ if $(p^k-1)|n(p^j-1)$.  We call elements satisfying
(\ref{eq:homogen}) ``constants.''

Constants behave as one might expect.  Adding a constant obviously does not change the derivative; multiplying by a constant $m$ scales the derivative by the same amount:
\alx{
    D_q (mf) &= f D_q(m) + M_q(m)D_q(f)  \\
             &= f\cdot 0  + mD_q(f)  \\
             &= m D_q(f)
}

\section{The exponential function exp}

Consider the equation $D_q x^e = x^e$.  Then
    $$D_q x^e = [e]x^{e-1} = x^e$$ $$[e]x^{-1}=x^0$$
    $$\frac{x^{es}-x^0}{x^s-x}=x^0$$ \al{ x^{es} &= x^s-x+1
    \label{eq:exp} } so if $\theta^s-\theta+1$ is generated by
$\theta^s$ then the equation will hold for at least one $e$.  We may
then define the function $\exp = x^e$; there is no reason to prefer one
solution over another.

We use $\exp$ to illustrate a subtlety of the chain rule.  One might conclude that $D_q (\exp^m)=[m]x^{me+m-1}$:
\al{
    D_q x^{me} &= D_q (x^e(x^m)) \nonumber \\
               &= D_q (x^e)(x^m) \cdot D_q(x^m) \nonumber\\
               &= x^e(x^m) \cdot [m]x^{m-1} \label{eq:bad} \\
               &= x^{me} \cdot [m]x^{m-1} \nonumber \\
               &= [m]x^{me+m-1} \nonumber
}
but (\ref{eq:bad}) does not follow.  It is only when applied directly
to $\theta$ that $D_qx^e=x^e$.  Here, $D_qx^e$ is applied to the
function $x^m$ and then to $\theta$.

The true equation may be found by examining the derivatives of the first few powers of $\exp$:
\alx{
    D_q (\exp^2) &= D_q (x^{2e}) \\
               &= D_q (x^e \cdot x^e) \\
               &= x^e D_q(x^e) + M_q(x^e)D_q(x^e) \\
               &= x^{2e} + q^e x^{2e} \\
               &= (q^0(x^e)+q^1(x^e))\cdot (x^e)^2 \\
               &= ([2]x^2)(x^e)
}\alx{
    D_q (\exp^3) &= D_q (x^{3e}) \\
               &= D_q (x^e\cdot x^{2e}) \\
               &= x^{2e} D_q(x^e)+ M_q(x^e)D_q(x^{2e})  \\
               &= x^{3e} + q^e\cdot (q^0+q^e) x^{3e} \\
               &= (q^0(x^e) + q^1(x^e) + q^2(x^e)) (x^e)^3 \\
               &= ([3]x^3)(x^e)
}
The pattern is immedately clear: $D_q(\exp^m)= ([m]x^m)(\exp)$, as one
would hope.

We can now prove the result by induction.  Assume that $D_q(\exp^{(m-1)})$ is of the form $\nobreak{([m-1]x^{m-1})(\exp)}$.  Then
\alx{
    D_q (\exp^m) &= D_q(x^{me}) \\
               &= D_q(x^ex^{(m-1)e}) \\
               &= x^{(m-1)e} D_q(x^e) + M_q(x^e)D_q(x^{(m-1)e}) \\
               &= x^{me} + q^ex^e \cdot ([m-1]x^{m-1})(x^e) \\
               &= ((q^0+q[m-1])x^m)(x^e)\\
               &= ([m]x^m)(x^e) \\
               &= ([m]x^m)(\exp).
}

\section{Commutation}

As with the standard $q$-derivative, $[D_q,x\cdot] = M_q$:
\alx{
    [D_q,x\cdot](f)  &= D_q(xf) - xD_q(f) \\
              &= fD_q(x) + M_q(x)D_q(f) - xD_q(f) \\
              &= f + qxD_q(f) - xD_q(f) \\
              &= f + d_q(x)D_q(f) \\
              &= (x + d_q(x)D_q)(f) \\
              &= (x + d_q(x)\frac{d_q}{d_q(x)})(f) \\
              &= (x + d_q)(f) \\
              &= M_q(f)
}

If we define the $q$-commutator $[f,g]_q \equiv fg-M_q(g)f$, then we find that 
\alx{
    [D_q,x\cdot]_q(f) &= D_q(xf) - M_q(x)D_q(f) \\
                &= fD_q(x) + M_q(x)D_q(f) - M_q(x)D_q(f)\\
                &= f.
}

We can define a Hamiltonian operator via the anticommutator $H=\lbrace D_q, x\cdot \rbrace$ to get
\alx{
    Hf &= D_q(xf) + xD_q(f) \\
       &= fD_q(x) + qxD_q(f) + xD_q(f) \\
       &= f + [2]xD_q(f),
}
so the ``energy'' of a finite field element $x^n(\theta)$ is
\alx{
    Hx^n &= x^n+[2]xD_q(x^n) \\
         &= (1+[2][n]) x^n \\
         &= (1+(1+q)[n]) x^n \\
         &= (1+[n]+q[n]) x^n \\
         &= ([n+1]+[n]) x^n
}

\section{$q$-Antiderivative}

The $q$-derivative of a finite field element is an element itself.  If
we add the constant $1$, the derivative does not change, so at most
half of the elements have antiderivatives.  If an element has an
antiderivative, then it is unique up to an additive constant: suppose
$f$ has two antiderivatives $F_1$ and $F_2$.  Then let $\phi=F_1-F_2$.
Now $D_q(\phi)=0$; but any function for which that holds true is a
constant by definition.

The integral operator $\int_q(d_q\cdot)$ is the Moore--Penrose inverse
of $D_q$.  Thus the equation $D_q(F)=f$ has a solution iff $f =
D_q(\int_q(fd_q)))$.

\section{Higher derivatives}
Because $[n]$ is a function of $x$, there are correction terms on the higher derivatives.  For instance,
\alx{
    D_q^2(x^n) &= D_q(D_q(x^n)) \\
               &= D_q([n]x^{n-1})\\
               &= [n]D_q(x^{n-1}) + M_q(x^{n-1})D_q([n]) \\
               &= [n][n-1]x^{n-2} + (qx)^{n-1}D_q([n])
}

It is these extra terms that give rise to trigonometric-like
functions.  We've already seen exp; there are others like sinh and cosh
with larger periods.

There will be a subspace, however, for which iterated derivatives
eventually yield zero.  This subspace always includes the vectors
$\lbrace x, 1\rbrace$, and may include more.

We can define an inner product in this subspace.  Let $J_q =
\int_q(d_q\cdot)$ and without loss of generality, let $n\geq m$.  Then
    $$\langle J_q^n,J_q^m\rangle = \langle 1,D_q^nJ_q^m\rangle =
    \delta_{n,m}.$$ The function exp is an eigenvector of $D_q$, so it
is orthogonal to the subspace:
    $$\langle J_q^n, \exp\rangle = \langle 1, D_q^n \exp\rangle =
    \langle 1,\exp\rangle = 0.$$ Other trigonometric functions are
    defined by the period with which they repeat.  sinh, for example,
is an eigenvector of $D_q^2$, and a similar identity holds.

\section{Logarithmic $q$-derivative}
The logarithmic derivative is defined as
    $$ \ld \equiv \frac{D_q}{x}. $$
The logarithmic derivative of a product of terms is the $q$-deformed sum of the logarithmic derivatives of the terms:
\alx{
    \ld(x^n \cdot x^m) &= \frac{D_q(x^n \cdot x^m)}{x^n\cdot x^m} \\
                       &= \frac{M_q(x^n) D_q(x^m) + x^m D_q (x^n)}{x^n\cdot x^m}\\
                       &= \frac{M_q(x^n)}{x^n}\ld(x^m) + \ld(x^n) \\
                       &= q^n \ld(x^m) + \ld(x^n) \\
                       &= \ld(x^m) + q^m \ld(x^n)
}

The logarithmic derivative of powers of exp has a nice form:
    $$ \ld(\exp^n) = \frac{D_q(\exp^n)}{\exp^n} =
    \frac{([n]x^n)(\exp)}{\exp^n} = [n](\exp) $$ which suggests a
``natural discrete $q$-logarithm'' for finite field elements.  However,
while this logarithm is easy to compute, the $q$-deformed
multiplication necessary to solve the Diffie-Hellman problem is hard.

\section{Examples}
We consider the field GF($2^4$) with the field polynomial $x^4-x-1$.
There are three possible values $q$ may take: $x^1, x^3,$ and $x^7$.
Each gives rise to different structures.

\ctr{\begin{tabular}{lllll}
    $n$  &$\theta^n$ & $D_x(\theta^n)$ & $D_{x^3}(\theta^n)$ & $D_{x^7}(\theta^n)$\\
    \hline \\
    $-\oo$ & 0000  & 0000 & 0000 & 0000 \\
    0      & 0001  & 0000 & 0000 & 0000 \\
    1      & 0010  & 0001 & 0001 & 0001 \\
    4      & 0011  & 0001 & 0001 & 0001 \\
    2      & 0100  & 0110 & 0001 & 0111 \\
    8      & 0101  & 0110 & 0001 & 0111 \\
    5      & 0110  & 0111 & 0000 & 0110 \\
    10     & 0111  & 0111 & 0000 & 0110 \\
    3      & 1000  & 1111 & 0111 & 0101 \\
    14     & 1001  & 1111 & 0111 & 0101 \\
    9      & 1010  & 1110 & 0110 & 0100 \\
    7      & 1011  & 1110 & 0110 & 0100 \\
    6      & 1100  & 1001 & 0110 & 0010 \\
    12     & 1101  & 1001 & 0110 & 0010 \\
    11     & 1110  & 1000 & 0111 & 0011 \\
    13     & 1111  & 1000 & 0111 & 0011
\end{tabular}}

\subsection{$q=x$} 
We have constants 0, 1.  ``Trig'' functions include
$\theta^{10}=0111=\exp$, $\theta^{13}=1111=\sinh$, and
$\theta^3=1000=\cosh$.  The names we've chosen are fairly arbitrary;
they are only meant to reflect the period with which the derivative
returns to itself.  The element $\theta$ has no antiderivative, so we
have an inner product acting on the subspace $\lbrace 1, x \rbrace$ of
the space $\lbrace 1, x, \exp, \sinh \rbrace$.

\subsection{$q=x^3$} Nonzero constants are $\theta^0=1, \theta^5=0110,$
and $\theta^{10}=0111$, the cube roots of 1.  There are no trig
functions.  A basis for the space is $\lbrace 1, x \rbrace$.

\subsection{$q=x^7$} We have the constants 0, 1 and the trig function
exp.  In this case, $J_q\theta=\theta^6$, so we have the inner product
on a three-dimensional subspace $\lbrace 1, x, x^6 \rbrace$, while the
complete basis is $\lbrace 1, x, x^6, \exp \rbrace$.

\bigskip
\hrule
\bigskip

\noindent 2000 {\it Mathematics Subject Classification}:
Primary 05A30; Secondary 11T99 .

\noindent \emph{Keywords:} $q$-calculus,
number derivative, arithmetic derivative. 

\bigskip
\hrule
\bigskip

\vspace*{+.1in}
\noindent
Received August 24 2004;
revised version received January 13 2005.  
Published in {\it Journal of Integer Sequences},
January 14 2005.
\bigskip
\hrule
\bigskip

\noindent
Return to
\htmladdnormallink{Journal of Integer Sequences home page}{http://www.math.uwaterloo.ca/JIS/}.
\vskip .1in

\end{document}